\documentclass[11pt]{amsart}
\usepackage{macros}

\title{On a theorem of Baxter and Zeilberger via a result of Roselle}
\author{Joshua P. Swanson}
\date{\today}

\begin{document}

\begin{abstract}
  We provide a new proof of a result of Baxter and Zeilberger showing that $\inv$ and
  $\maj$ on permutations are jointly independently asymptotically normally distributed.
  The main feature of our argument is that it uses a generating function due to
  Roselle, answering a question raised by Romik and Zeilberger.
\end{abstract}
\maketitle

\section{Introduction}

For a permutation $w = w_1 \cdots w_n \in S_n$, the \textit{inversion} and
\textit{major index} statistics are given by
  \[ \inv(w) := \#\{i < j : w_i > w_j\} \qquad\text{and}\qquad
      \maj(w) := \sum_{\substack{1 \leq i \in \leq n-1\\w_i > w_{i+1}}} i. \]
It is well-known that $\inv$ and $\maj$ are equidistributed on $S_n$ with common mean
and standard deviation
  \[ \mu_n = \frac{n(n-1)}{4} \qquad\text{and}\qquad
      \sigma_n^2 = \frac{2n^3 + 3n^2 - 5n}{72}. \]
(These results also follow easily from our arguments.) In \cite{1004.1160}, Baxter and
Zeilberger proved that $\inv$ and $\maj$ are jointly
independently asymptotically normally distributed as $n \to \infty$. More precisely,
define normalized random variables on $S_n$
\begin{equation}\label{eq:XnYn}
  X_n := \frac{\inv - \mu_n}{\sigma_n}, \qquad Y_n := \frac{\maj - \mu_n}{\sigma_n}.
\end{equation}

\begin{Theorem}[Baxter--Zeilberger, \cite{1004.1160}]\label{thm:BZ}
  For each $u, v \in \bR$, we have
    \[ \lim_{n \to \infty} \bP[X_n \leq u, Y_n \leq v]
        = \frac{1}{2\pi} \int_{-\infty}^u \int_{-\infty}^v e^{-x^2/2} e^{-y^2/2}\,dy\,dx. \]
\end{Theorem}

See \cite{1004.1160} for further historical background. Baxter and Zeilberger's argument
involves mixed moments and recurrences based on combinatorial manipulations with
permutations. Romik suggested a generating function due to Roselle, quoted as
\Cref{thm:Roselle} below, should provide another approach. Zeilberger subsequently
offered a \$300 reward for such an argument. The aim of this note is to give such
a proof. Our overarching motivation is to give a \textit{local limit theorem},
i.e.~a formula for the counts $\#\{w \in S_n : \inv(w) = u, \maj(w) = v\}$,
with an explicit error term, which will be the subject of a future article.
For further context, see \cite{dz} and \cite{MR3482667}.

\section{Consequences of Roselle's Formula}

Here we recall Roselle's formula, originally stated in different but equivalent terms, and
derive a generating function expression which quickly motivates \Cref{thm:BZ}.

\begin{Definition}
  Let $H_n$ be the bivariate $\inv, \maj$ generating function on $S_n$, i.e.
    \[ H_n(p, q) := \sum_{w \in S_n} p^{\inv(w)} q^{\maj(w)}. \]
\end{Definition}

\begin{Theorem}[Roselle, \cite{MR0342406}]\label{thm:Roselle}
  We have
  \begin{equation}\label{eq:roselle}
    \sum_{n \geq 0} \frac{H_n(p, q) z^n}{(p)_n (q)_n}
        = \prod_{a, b \geq 0} \frac{1}{1 - p^a q^b z}
  \end{equation}
  where $(p)_n := (1-p)(1-p^2) \cdots (1-p^n)$.
\end{Theorem}

The following is the main result of this section.

\begin{Theorem}\label{thm:Hncla}
  There are constants $c_\mu \in \bZ$ indexed by integer partitions $\mu$ such that
  \begin{equation}\label{eq:HnFn}
    \frac{H_n(p, q)}{n!} = \frac{[n]_p! [n]_q!}{n!^2} F_n(p, q)
  \end{equation}
  where
  \begin{equation}\label{eq:Fncmu}
    F_n(p, q) = \sum_{d=0}^n [(1-p)(1-q)]^d
      \sum_{\substack{\mu \vdash n\\\ell(\mu) = n-d}}
      \frac{c_\mu}{\prod_i [\mu_i]_p [\mu_i]_q}
  \end{equation}
  and $[n]_p! := [n]_p [n-1]_p \cdots [1]_p$, $[c]_p := 1 + p + \cdots + p^{c-1}
  = (1-p^c)/(1-p)$.
\end{Theorem}

An explicit expression for $c_\mu$ is given below in \eqref{eq:cmu}. The rest of this section is
devoted to proving \Cref{thm:Hncla}. Straightforward manipulations with \eqref{eq:roselle}
immediately yield \eqref{eq:HnFn} where
\begin{equation}\label{eq:FnRatio}
   F_n(p, q) := (1-p)^n (1-q)^n n! \cdot
     \{z^n\} \left(\prod_{a, b \geq 0} \frac{1}{1 - p^a q^b z}\right)
\end{equation}
and $\{z^n\}$ here refers to extracting the coefficient of $z^n$.
Thus it suffices to show \eqref{eq:FnRatio} implies \eqref{eq:Fncmu}.
By standard arguments, the $z^n$ coefficient of the product over $a, b$ in \eqref{eq:FnRatio}
is the bivariate generating function of size-$n$ multisets of pairs $(a, b) \in \bZ_{\geq 0}^2$,
where the weight of such a multset is its sum.

\begin{Definition}
  For $\lambda \vdash n$, let $M_\lambda$ be the bivariate generating function for
  multisets of pairs $(a, b) \in \bZ_{\geq 0}^n$ of type $\lambda$, i.e.~some element has
  multiplicity $\lambda_1$, another element has multiplicity $\lambda_2$, etc.
\end{Definition}

We clearly have
\begin{equation}\label{eq:FnMla}
  \{z^n\}\left(\prod_{a, b \geq 0} \frac{1}{1 - p^a q^b z}\right)
  = \sum_{\lambda \vdash n} M_\lambda(p, q),
\end{equation}
though the $M_\lambda$ are inconvenient to work with, so we perform a change of basis.

\begin{Definition}
  Let $P[n]$ denote the lattice of set partitions of $[n] := \{1, 2, \ldots, n\}$ with
  minimum $\widehat{0} = \{\{1\}, \{2\}, \ldots, \{n\}\}$ and maximum
  $\widehat{1} = \{\{1, 2, \ldots, n\}\}$. Here $\Lambda \leq \Pi$ means that
  $\Pi$ can be obtained from $\Lambda$ by merging blocks of $\Lambda$. The
  \textit{type} of a set partition $\Lambda$ is the integer partition
  obtained by rearranging the list of the block sizes of $\Lambda$ in weakly decreasing order.
  For $\lambda \vdash n$, set
    \[ \Lambda(\lambda) := \{\{1, 2, \ldots, \lambda_1\},
        \{\lambda_1+1, \lambda_1+2, \ldots, \lambda_1+\lambda_2\}, \ldots\}, \]
  which has type $\lambda$.
\end{Definition}

\begin{Definition}
  For $\Pi \in P[n]$, let $R_\Pi$ denote the bivariate generating function for lists $L
  \in (\bZ_{\geq 0}^2)^n$ where for each block of $\Pi$ the entries in $L$ from that block
  are all equal. Similarly, let $S_\Pi$ denote the bivariate generating function of lists $L$ where
  in addition to entries from the same block being equal, entries from two different blocks
  are not equal.
\end{Definition}

We easily see that
\begin{equation}\label{eq:Rla}
  R_\Lambda(p, q) = \prod_{A \in \Lambda} \frac{1}{(1 - p^{\#A})(1 - q^{\#A})}
\end{equation}
and that
\begin{equation}\label{eq:RlaSpi}
  R_\Lambda(p, q) = \sum_{\Pi : \Lambda \leq \Pi} S_\Pi,
\end{equation}
so that, by M\"obius inversion on $P[n]$,
\begin{equation}\label{eq:SpiRla}
  S_\Pi = \sum_{\Lambda : \Pi \leq \Lambda} \mu(\Pi, \Lambda) R_\Lambda.
\end{equation}
Under the ``forgetful'' map from lists to multisets, a multiset of type $\lambda \vdash n$ has
fiber of size $\binom{n}{\lambda}$. It follows that
\begin{equation}\label{eq:SpiMla}
  S_{\Pi(\lambda)} = \frac{n!}{\lambda!} M_\lambda
\end{equation}
where $\lambda! := \lambda_1! \lambda_2! \cdots$. Combining in order \eqref{eq:FnRatio},
\eqref{eq:FnMla}, \eqref{eq:SpiMla}, \eqref{eq:SpiRla}, and \eqref{eq:Rla} gives
\begin{equation}\label{eq:FnRla}
  F_n(p, q) = \sum_{d=0}^n [(1-p)(1-q)]^d \sum_{\lambda \vdash n} \lambda!
    \sum_{\substack{\Lambda : \Pi(\lambda) \leq \Lambda\\\#\Lambda = n-d}}
      \frac{\mu(\Pi(\lambda), \Lambda)}{\prod_{A \in \Lambda} [\#A]_p [\#A]_q}.
\end{equation}
Now \eqref{eq:Fncmu} follows from \eqref{eq:FnRla} where
\begin{equation}\label{eq:cmu}
  c_\mu = \sum_{\lambda \vdash n} \lambda!
    \sum_{\substack{\Lambda : \Pi(\lambda) \leq \Lambda\\\type(\Lambda)=\mu}}
    \mu(\Pi(\lambda), \Lambda).
\end{equation}
This completes the proof of \Cref{thm:Hncla}.

\begin{Remark}\label{rem:macmahon}
  From \eqref{eq:cmu}, $c_{(1^n)} = 1$ since the sum only involves
  $\Lambda = \widehat{0}$. Letting $p \to 1$ in \eqref{eq:Fncmu},
  the only surviving term is $d=0$ and $\lambda = (1^n)$. Consequently,
  $H_n(1, q) = [n]_q!$, recovering a classic result of MacMahon
  \cite[\S1]{MR1576566}.
\end{Remark}

\begin{Remark}
  Using \eqref{eq:HnFn}, we see that the probability generating function
  (discussed below in \Cref{ex:pgfchar}) $H_n(p, q)/n!$ differs from
  $[n]_p! [n]_q!/n!^2$ by precisely the correction factor $F_n(p, q)$. Using
  \eqref{eq:FnRatio}, this factor has the following combinatorial interpretation:
    \[ F_n = \frac{n! \cdot \text{g.f. of size-$n$ multisets from $\bZ_{\geq 0}^2$}}
                                 {\text{g.f. of size-$n$ lists from $\bZ_{\geq 0}^2$}}. \]
  Intuitively, the numerator and denominator should be the same ``up to first order.''
  \Cref{thm:FnBound} will give one precise sense in which they are asymptotically equal.
\end{Remark}

\section{Estimating the Correction Factor}

This section is devoted to showing that the correction factor $F_n(p, q)$ from
\Cref{thm:Hncla} is negligible in an appropriate sense, \Cref{thm:FnBound}. Recall
that $\sigma_n$ denotes the standard deviation of $\inv$ or $\maj$ on $S_n$.

\begin{Theorem}\label{thm:FnBound}
  Uniformly on compact subsets of $\bR^2$, we have
    \[ F_n(e^{is/\sigma_n}, e^{it/\sigma_n}) \to 1
        \qquad \text{as} \qquad n \to \infty \]
\end{Theorem}

We begin with some simple estimates starting from \eqref{eq:FnRla}
which motivate the rest of the inequalities in this section.
We may assume $|s|, |t| \leq M$ for some fixed $M$. Setting
$p=e^{is/\sigma_n}, q=e^{it/\sigma_n}$, we have
$|1-p| = |1-\exp(is/\sigma_n)| \leq |s|/\sigma_n$. For $n$ sufficiently
large compared to $M$, we also have $|s/\sigma_n| \ll 1$ and so, for all
$c \in \bZ_{\geq 1}$, $|[c]_p| = |[c]_{\exp(is/\sigma_n)}| \geq 1$.
Thus for $n$ sufficiently large, \eqref{eq:FnRla} gives
\begin{equation}\label{eq:ModFnRla}
  |F_n(e^{is/\sigma_n}, e^{it/\sigma_n}) - 1|
    \leq \sum_{d=1}^n \frac{|st|^d}{\sigma_n^{2d}}
      \sum_{\lambda \vdash n} \lambda!
      \sum_{\substack{\Lambda : \Pi(\lambda) \leq \Lambda\\\#\Lambda = n-d}}
        |\mu(\Pi(\lambda), \Lambda)|
\end{equation}

\begin{Lemma}
  Suppose $\lambda \vdash n$ with $\ell(\lambda) = n-k$, and fix $d$. Then
  \begin{equation}\label{eq:dksum}
    \sum_{\substack{\Lambda : \Pi(\lambda) \leq \Lambda\\\#\Lambda = n-d}}
      \mu(\Pi(\lambda), \Lambda)
      = (-1)^{d-k} \sum_{\substack{\Lambda \in P[n-k]\\\#\Lambda = n-d}}
         \prod_{A \in \Lambda} (\#A-1)!
  \end{equation}
  and the terms on the left all have the same sign $(-1)^{d-k}$. The sums are
  empty unless $n \geq d \geq k \geq 0$.
  
  \begin{proof}
    The upper order ideal $\{\Lambda \in P[n] : \Pi(\lambda) \leq \Lambda\}$ is isomorphic
    to $P[n-k]$ by collapsing the $n-k$ blocks of $\Pi(\lambda)$ to singletons. This
    isomorphism preserves the number of blocks. Furthermore, recall that in $P[n]$ we
    have
      \[ \mu(\widehat{0}, \widehat{1}) = (-1)^{n-1} (n-1)!, \]
    from which it follows easily that
    \begin{equation}\label{eq:mu0La}
      \mu(\widehat{0}, \Lambda) = \prod_{A \in \Lambda} (-1)^{\#A - 1} (\#A - 1)!.
    \end{equation}
    The result follows immediately upon combining these observations.
  \end{proof}
\end{Lemma}

\begin{Lemma}
  Let $\lambda \vdash n$ with $\ell(\lambda) = n-k$ and $n \geq d \geq k \geq 0$. Then
  \begin{equation}\label{eq:ladk_bound}
    \sum_{\substack{\Lambda : \Pi(\lambda) \leq \Lambda\\\#\Lambda = n-d}}
      |\mu(\Pi(\lambda), \Lambda)| \leq (n-k)^{2(d-k)}.
  \end{equation}
  
  \begin{proof}
    Using \eqref{eq:dksum}, we can interpret the sum as the number of permutations
    of $[n-k]$ with $n-d$ cycles, which is a Stirling number of the first kind. There are
    well-known asymptotics for these numbers, though the stated elementary bound suffices for
    our purposes. We induct on $d$. At $d=k$, the result is trivial.
    Given a permutation of $[n-k]$ with $n-d$ cycles, choose $i, j \in [n-k]$ from different
    cycles. Suppose the cycles are of the form $(i'\ \cdots\ i)$ and $(j\ \cdots\ j')$.
    Splice the two cycles together to obtain
      \[ (i'\ \cdots\ i\ j\ \cdots\ j'). \]
    This procedure constructs every permutation of $[n-k]$ with $n-(d+1)$ cycles and
    requires no more than $(n-k)^2$ choices. The result follows.
  \end{proof}
\end{Lemma}

\begin{Lemma}
  For $n \geq d \geq k \geq 0$, we have
  \begin{equation}\label{eq:dksum_bound}
    \sum_{\substack{\lambda \vdash n\\\ell(\lambda) = n-k}}
       \lambda! \sum_{\substack{\Lambda : \Pi(\lambda) \leq \Lambda\\\#\Lambda = n-d}}
       |\mu(\Pi(\lambda), \Lambda)|
       \leq (n-k)^{2d-k} (k+1)!.
  \end{equation}
  
  \begin{proof}
      For $\lambda \vdash n$ with $\ell(\lambda) = n-k$, $\lambda!$ can be thought of as the
      product of terms obtained from filling the $i$th row of $\lambda$ with $1, 2, \ldots,
      \lambda_i$. Alternatively, we may fill the cells of $\lambda$ as follows: put $n-k$ one's in
      the first column, and fill the remaining cells with the numbers $2, 3, \ldots, k+1$
      starting at the
      largest row and proceeding left to right. It's easy to see the labels of the first filling
      are bounded above by the labels of the second filling, so that $\lambda! \leq (k+1)!$.
      Furthermore, each $\lambda \vdash n$ with $\ell(\lambda) = n-k$ can be constructed
      by first placing $n-k$ cells in the first column and then deciding on which of the $n-k$ rows
      to place each of the remaining $k$ cells, so there are no more than $(n-k)^k$ such
      $\lambda$. The result follows from combining these bounds with
      \eqref{eq:ladk_bound}.
  \end{proof}
\end{Lemma}

\begin{Lemma}\label{lem:d_bound}
  For $n$ sufficiently large, for all $0 \leq d \leq n$ we have
  \begin{align*}
    \sum_{\lambda \vdash n} \lambda!
      \sum_{\substack{\Lambda : \Pi(\lambda) \leq \Lambda\\\#\Lambda = n-d}}
      |\mu(\Pi(\lambda), \Lambda)| \leq 3n^{2d}.
  \end{align*}
  
  \begin{proof}
    For $n \geq 2$ large enough, for all $n \geq k \geq 2$ we see that $(k+1)! < n^{k-1}$.
    Using \eqref{eq:dksum_bound} gives
    \begin{align*}
      \sum_{\lambda \vdash n} \lambda!
      \sum_{\substack{\Lambda : \Pi(\lambda) \leq \Lambda\\\#\Lambda = n-d}}
      |\mu(\Pi(\lambda), \Lambda)|
        &\leq \sum_{k=0}^d (n-k)^{2d-k} (k+1)! \\
        &\leq n^{2d} + 2(n-1)^{2d-1} + \sum_{k=2}^d (n-k)^{2d-k} n^{k-1} \\
        &\leq n^{2d} + 2n^{2d-1} + \sum_{k=2}^d n^{2d-1} \\
        &= n^{2d} + 2n^{2d-1} + (d-1) n^{2d-1}
        \leq 3n^{2d}.
    \end{align*}
  \end{proof}
\end{Lemma}

We may now complete the proof of \Cref{thm:FnBound}. Combining \Cref{lem:d_bound}
and \eqref{eq:ModFnRla} gives
  \[ |F_n(e^{is/\sigma_n}, e^{it/\sigma_n}) - 1| \leq 3\sum_{d=1}^n
      \frac{(Mn)^{2d}}{\sigma_n^{2d}}. \]
Since $\sigma_n^2 \sim n^3/36$ and $M$ is constant,
$(Mn)^{2d}/\sigma_n^{2d} \sim (36^2M^2/n)^d$. Since $M$ is constant,
using a geometric series it follows that
  \[ \lim_{n \to \infty} \sum_{d=1}^n \frac{(Mn)^{2d}}{\sigma_n^{2d}} = 0, \]
completing the proof of \Cref{thm:FnBound}.

\begin{Remark}
  Indeed, the argument shows that
  $|F_n(e^{is/\sigma_n}, e^{it/\sigma_n}) - 1| = O(1/n)$. The above estimates
  are particularly far from sharp for large $d$, though for small $d$ they are
  quite accurate. Working directly with \eqref{eq:FnRla}, one finds the $d=1$
  contribution to be
  \begin{align*}
    (1-p)(1-q)
    \frac{2 - \binom{n}{2}}{[2]_p [2]_q}.
  \end{align*}
  Letting $p = e^{is/\sigma_n}, q = e^{it/\sigma_n}$, straightforward estimates
  shows that this is $\Omega(1/n)$. Consequently, the preceding arguments
  are strong enough to identify the leading term, and in particular
   \[ |F_n(e^{is/\sigma_n}, e^{it/\sigma_n}) - 1| = \Theta(1/n). \]
\end{Remark}

\section{Deducing Baxter and Zeilberger's Result}\label{sec:cfs}

We next summarize enough of the standard theory of characteristic functions to prove
\Cref{thm:BZ} using \eqref{eq:HnFn} and \Cref{thm:FnBound}.

\begin{Definition}
  The \textit{characteristic function} of an $\bR^k$-valued random variable
  $X = (X_1, \ldots, X_k)$ is the function $\phi_X \colon \bR^k \to \bC$ given by
    \[ \phi_X(s_1, \ldots, s_k) := \bE[\exp(i(s_1 X_1 + \cdots + s_k X_k))]. \]
\end{Definition}

\begin{Example}\label{ex:normchar}
  It is well-known that the characteristic function of the standard normal random variable with
  density $\frac{1}{\sqrt{2\pi}} e^{-x^2/2}$ is $e^{-s^2/2}$. Similarly, the characteristic
  function of a bivariate jointly independent standard normal random variable with density
  $\frac{1}{2\pi} e^{-x^2/2 - y^2/2}$ is $e^{-s^2/2 - t^2/2}$.
\end{Example}

\begin{Example}\label{ex:pgfchar}
  If $W$ is a finite set and $\stat = (\stat_1, \ldots, \stat_k) \colon W \to \bZ_{\geq 0}^k$ is some statistic,
  the \textit{probability generating function} of $\stat$ on $W$ is
    \[ P(x_1, \ldots, x_k) := \frac{1}{\#W}
        \sum_{w \in W} x_1^{\stat_1(w)} \cdots x_k^{\stat_k(w)}. \]
  The characteristic function of the corresponding random variable $X$
  where the $w$ are chosen uniformly from $W$ is
    \[ \phi_X(s_1, \ldots, s_k) = P(e^{is_1}, \ldots, e^{is_k}). \]
\end{Example}

From \Cref{ex:pgfchar}, \Cref{rem:macmahon}, and an easy calculation, it follows that the
characteristic functions of the random variables $X_n$ and $Y_n$ from \eqref{eq:XnYn} are
\begin{equation}\label{eq:XnYnChar}
  \phi_{X_n}(s) = e^{-i\mu_n s/\sigma_n} \frac{[n]_{e^{is/\sigma_n}}!}{n!}
    \qquad\text{and}\qquad
  \phi_{Y_n}(t) = e^{-i\mu_n t/\sigma_n} \frac{[n]_{e^{it/\sigma_n}}!}{n!}
\end{equation}

An analogous calculation for the random variable $(X_n, Y_n)$ together with
\eqref{eq:XnYnChar} and \eqref{eq:HnFn} gives
\begin{equation}\label{eq:XnYnChar2}
  \begin{split}
    \phi_{(X_n, Y_n)}(s, t) &= e^{-i(\mu_n s/\sigma_n + \mu_n t/\sigma_n)}
      \frac{H_n(e^{is/\sigma_n}, e^{it/\sigma_n})}{n!} \\
    &= \phi_{X_n}(s) \phi_{Y_n}(t) F_n(e^{is/\sigma_n}, e^{it/\sigma_n}).
  \end{split}
\end{equation}

\begin{Theorem}[Multivariate L\'evy Continuity, {\cite[p.~383]{MR1324786}}]\label{thm:levy}
  Suppose that $X^{(1)}$, $X^{(2)}$, $\ldots$ is a sequence of $\bR^k$-valued
  random variables
  and $X$ is an $\bR^k$-valued random variable. Then $X^{(1)}, X^{(2)}, \ldots$ converges
  in distribution to $X$ if and only if $\phi_{X^{(n)}}$ converges pointwise to $\phi_X$.
\end{Theorem}

If the distribution function of $X$ is continuous everywhere, convergence in distribution means
that for all $u_1, \ldots, u_k$ we have
  \[ \lim_{n \to \infty} \bP[X^{(n)}_i \leq u_i, 1 \leq i \leq k]
      = \bP[X_i \leq u_i, 1 \leq i \leq k]. \]
Many techniques are available for proving that $\inv$ and $\maj$ on $S_n$ are asymptotically
normal. The result is typically attributed to Feller.

\begin{Theorem}{\cite[p.~257]{MR0228020}}\label{thm:feller}
  The sequences of random variables $X_n$ and $Y_n$ from \eqref{eq:XnYn} each converge
  in distribution to the standard normal random variable.
\end{Theorem}

We may now complete the proof of \Cref{thm:BZ}. From \Cref{thm:feller} and
\Cref{ex:normchar}, we have for all $s, t \in \bR$
\begin{equation}\label{eq:XnYnLim}
  \lim_{n \to \infty} \phi_{X_n}(s) = e^{-s^2/2}
  \qquad\text{and}\qquad
  \lim_{n \to \infty} \phi_{Y_n}(t) = e^{-t^2/2}.
\end{equation}
Combing in order \eqref{eq:XnYnLim}, \eqref{eq:XnYnChar2}, and \Cref{thm:FnBound}
gives
  \[ \lim_{n \to \infty} \phi_{(X_n, Y_n)}(s, t) = e^{-s^2/2 - t^2/2}. \]
\Cref{thm:BZ} now follows from \Cref{ex:normchar} and \Cref{thm:levy}.

\section{Acknowledgments}

The author would like to thank Dan Romik and Doron Zeilberger for
providing the impetus for the present work and feedback on the
manuscript. He would also like to thank Sara Billey and Matja\v{z}
Konvalinka for valuable discussion on related work, and he gratefully
acknowledges Sara Billey for her very careful reading of the manuscript and
many helpful suggestions.

\bibliography{zeilberger_300}{}
\bibliographystyle{alpha}

\end{document}